\documentclass[10pt]{article}
\usepackage{amssymb}  
\usepackage{colordvi,times,fancybox,latexsym}
\usepackage{graphics}
\usepackage{graphicx}

\newenvironment{proof}{\noindent {\em Proof}.\ }{\hspace
*{\fill}\medskip}
\newtheorem{theorem}{Theorem}

\newtheorem{lemma}[theorem]{Lemma}

\newtheorem{remark}[theorem]{Remark}

\begin{document}

\title{Further Results on Active Magnetic Bearing Control with Input Saturation%
\thanks{%
This work was supported in part by NSF grants DMS-0424011 and CMS-0447576.}}
\author{F. Mazenc\thanks{\textit{Mazenc} is with the Projet MERE INRIA-INRA; UMR
Analyse des Syst\`{e}mes et Biom\'{e}trie INRA; 2, pl. Viala; 34060
Montpellier, France; Frederic.Mazenc@ensam.inra.fr.},\ M.S. de Queiroz\thanks{%
\textit{de Queiroz} is with the Department of Mechanical
Engineering, Louisiana State University, Baton Rouge, LA 70803-6413,
dequeiroz@me.lsu.edu. (Corresponding author)},\ M. Malisoff\thanks{\textit{%
Malisoff} is with the Department of Mathematics, Louisiana State University,
Baton Rouge, LA 70803-4918, malisoff@lsu.edu.}, and F. Gao\thanks{\textit{Gao%
} is with the Department of Mechanical Engineering, Louisiana State
University, Baton Rouge, LA 70803-6413, fgao1@lsu.edu.}}
\date{}
\maketitle

\begin{abstract}
We study the low-bias stabilization of active magnetic bearings (AMBs)
subject to voltage saturation based on a recently proposed model for the AMB
switching mode of operation. Using a forwarding-like approach, we construct
a stabilizing controller of arbitrarily small amplitude and a
control-Lyapunov function for the AMB\ dynamics. We illustrate our
construction using a numerical example.
\end{abstract}


\section{Introduction\label{intro}}

Active magnetic bearings (AMBs) are being employed in a variety of rotating
machines (e.g., compressors, milling spindles, and flywheels) in place of
conventional mechanical bearings. Due to the non-contact nature of the
magnetic bearings and rotor, AMBs have the unique ability to suspend loads
with no friction, operate rotors at higher speeds, and operate in
applications where the use of lubricants is prohibited. Since AMBs can be
actively controlled, they offer other potential advantages over mechanical
bearings, viz., the elimination of vibration through active damping, the
adjustment of the stiffness of the suspended load, and the automatic
balancing of rotors. AMBs are normally operated with all electromagnets
active at all times. An alternative mode of operation is to activate only
one electromagnet along each direction at any given time. Although posing
more difficulty to the control design, this switching operation prevents
opposing electromagnets from producing counteracting forces, and thus helps
reduce power consumption.

Typically, an AMB is operated by introducing a sufficiently high, fixed
magnetic flux in each electromagnet, which is referred to as the \textit{%
bias flux}. The bias value is normally set to a fraction of the saturation
flux of the electromagnet. This procedure facilitates the design of the AMB
`control' flux, which is superimposed on the bias flux. Specifically, this
conservative practice allows the system to be modeled by a controllable
linear system, thus, enabling the use of linear control design techniques;
see for example \cite{fujita,matsumura,mohamed1}. Although the bias flux
facilitates the control synthesis, it increases electric power losses in the
AMB system, causing rotor heating and affecting the machine efficiency.
While lowering or eliminating the bias flux is desirable in order to
minimize power losses, it enhances the AMB\ system nonlinearities and may
lead to a control saturation\footnote{%
Saturation may arise due to limits on the amplitude of the output voltage of
the power amplifier driving the electromagnets.} or singularity. Due to
these conflicting objectives, the design of AMB\ controllers with reduced
power loss is a challenging problem.

In this paper, we consider the problem of low-bias\footnote{%
By \emph{low bias}, we mean an AMB\ control law where the bias can be set to
an arbitrarily small positive constant without affecting the system
stability.} control of AMBs operating in the switching mode with the
constraint that the input voltages are amplitude limited. This problem was
previously addressed in \cite{tsiotras1,TA02} using the nested saturation
design method of Teel \cite{teel1,teel2}. In addition, an optimal solution
was sought to the power-loss minimization problem under voltage saturation
in \cite{hu}. (A\ comprehensive literature review of low-, asymptotic-zero-,
and zero-bias AMB\ controllers \textit{without} voltage saturation can be
found in \cite{motee cdc03,tsiotras4}.) Here, we pursue a different approach
to stabilizing the AMB model of \cite{TA02}. Namely, we use a \textit{%
forwarding-like} method \cite{mazenc} to design a control law of arbitrarily
small amplitude that renders the AMB system globally asymptotically stable
(GAS) to the origin. The main differences between our approach and \cite
{TA02} are: (i) we allow the use of different saturation functions (both
`hard' and `soft'), thus providing more flexibility in the control
implementation, and (ii) we provide an explicit construction for a `global'
control Lyapunov function (CLF) for the system. This contrasts with the
result of \cite{TA02}, which relies on the standard `hard' saturation and
whose Lyapunov-like function is only an `asymptotic' CLF for the system,
i.e., a CLF only when the system operates in a certain region of the state
space.

The rest of this paper is organized as follows. In Section
\ref{model}, we introduce the AMB model and motivate our
stabilization problem. In Section \ref{sec3}, we prove a lemma that
constructs a stabilizing controller of arbitrarily small amplitude
and a CLF for a three-dimensional chain of integrators. In Section
\ref{sec4}, we use our lemma to construct a CLF and corresponding
stabilizing feedback for the AMB model. We provide a numerical
example in Section \ref{sec5}, and we close in Section \ref{sec6}
with a summary of our work.


\section{AMB Model and Problem Statement\label{model}}

The original nonlinear electromechanical model of the one
degree-of-freedom AMB system shown in Figure 1 below can be
subdivided into the mechanical subsystem dynamics, the magnetic
force equation, and the electrical subsystem dynamics. The
mechanical subsystem is governed by
\begin{equation}
m\ddot{y}=\sum_{i=1}^{2}F_{i}(\Phi _{i})\,,\,\,\,\,\,  \label{mech}
\end{equation}
where $m$ is the rotor mass, $y\in {\mathbb R}$ represents the
position of the rotor center, $\Phi _{i}\in {\mathbb R}$ is the
magnetic flux in the $i$th electromagnet, $F_{i}\left( \Phi
_{i}\right) \in {\mathbb R}$ denotes the force produced by the $i$th
electromagnet, given by \cite{embook}
\begin{equation}
F_{i}=\frac{\left( -1\right) ^{i+1}\Phi _{i}^{2}}{\mu _{0}A},\quad i=1,2,
\label{force}
\end{equation}
$\mu _{0}$ is the permeability of air, and $A$ is the cross-sectional area
of the electromagnet. The electrical subsystem is governed by the equations
\cite{embook}
\begin{equation}
N\dot{\Phi}_{i}+R_{i}I_{i}=v_{i},\quad i=1,2,  \label{elec}
\end{equation}
where $N$ denotes the number of coil turns in the electromagnet,
$R_{i}$ is the resistance of the $i$th electromagnet coil, $v_{i}\in
{\mathbb R}$ is the input control voltage of the $i$th
electromagnet, $I_{i}\in {\mathbb R}$ is the current in the $i$th
electromagnet which is related to the flux according to
\cite{embook}
\begin{equation}
I_{i}=\frac{2\left( g_{0}+\left( -1\right) ^{i}y\right) \Phi _{i}}{\mu _{0}AN%
}\,\,,\,\,i=1,2,  \label{current-flux}
\end{equation}
and $g_{0}$ is the nominal air gap.

In this paper, we use the new form proposed in \cite{TA02,tsiotras4} for the
AMB\ dynamics in the switching mode of operation. We only outline the model
derivation here, and refer the reader to \cite{tsiotras4} for the details
and justification. Consider that the flux is given by
\begin{equation}
\Phi _{i}=\Phi _{0}+\phi _{i},\quad i=1,2  \label{flux}
\end{equation}
where $\Phi _{0}>0$ is the constant \textit{bias} flux and $\phi _{i}$ is
the \textit{control} flux. Let the \textit{generalized} control flux be
defined as
\begin{equation}
\phi =\phi _{1}-\phi _{2},  \label{phi}
\end{equation}
and consider the voltage switching strategy
\begin{equation}
\begin{array}{l}
v_{1}=v,\quad v_{2}=0\quad \,\,\,\,\,\mathrm{if\quad }\phi \geq 0 \\
v_{1}=0,\quad v_{2}=-v\quad \mathrm{if\quad }\phi <0
\end{array}
\label{switch}
\end{equation}
where $v$ is the \textit{generalized} control voltage. Based on (\ref{flux}%
)-(\ref{switch}), the AMB\ model (\ref{mech})-(\ref{elec}) has the
equivalent form
\begin{equation}
\ddot{y}=\frac{1}{m\mu _{0}A}\left( 2\bar{\Phi}_{0}\phi +\phi \left|
\phi \right| \right) ,\;\;\dot{\phi}=\frac{v}{N},  \label{mech
elec2}
\end{equation}
where $\bar{\Phi}_{0}=\Phi _{0}+\min \left\{ \phi _{1}(0),\phi
_{2}(0)\right\} $ and the coil resistance was neglected for simplicity.

Now, consider that the input voltages to the original AMB model are
amplitude limited, i.e., $\left| v_{i}\right| \leq v_{\max }$, $i=1,2$ where
$v_{\max }$ is the known limit. Defining the states $x_{1}=y$, $x_{2}=\dot{y}
$, and $x_{3}=\phi $ and the change of input $v=Nu$, the above model becomes
\begin{equation}
\dot{x}_{1}=x_{2},\;\;\dot{x}_{2}=\beta _{0}x_{3}+\beta _{1}x_{3}\left|
x_{3}\right| ,\;\;\dot{x}_{3}=u,  \label{xdot amb}
\end{equation}
where $\beta _{0}=2\bar{\Phi}_{0}/(m\mu _{0}A)$ and $\beta _{1}=1/(m\mu
_{0}A)$, with the input constraint $\left| u\right| \leq v_{\max }/N$.

Our goal in this work is to design a feedback control law $u(x)$, satisfying
the above saturation constraint, such that $x=0$ is GAS where $%
x=(x_{1},x_{2},x_{3})^{\top }$. The GAS property means there is a
continuously differentiable ($C^{1}$) function $V:{\mathbb R}^{3}\to
[0,\infty )$ that is radially unbounded and zero only at the origin,
and for which the derivative $\dot{V}$ along all trajectories of the
system in closed-loop with the controller $u(x)$ is negative
definite. A function $V$ satisfying these requirements for some
feedback $u(x)$ is called a CLF for (\ref{xdot amb}) \cite{krstic}.
See \cite{tsiotras4} for a discussion on the relation between the
global asymptotic stability of $x$ and the stability of the original
AMB states.

\section{Preliminary Result \label{sec3}}

In the following lemma, we use a forwarding-like approach based on \cite
{mazenc} to construct a CLF and corresponding stabilizing state feedback of
arbitrarily small amplitude for a three-dimensional chain of integrators. We
let $\sigma $ denote the standard saturation projecting ${\mathbb R}$ onto $%
[-1,+1]$, i.e., $\sigma (s)=s$ when $|s|\le 1$ and $\sigma (s)=\mathrm{sign}%
(s)$ otherwise. However, see Remark 1 below for a version of our
construction where the hard saturation $\sigma $ is replaced by a smoother
function.

\begin{lemma}
\rm%
\label{mainlemma} Consider the chain of integrators
\begin{equation}
\dot{\xi}_{1}=\xi _{2},\;\;\dot{\xi}_{2}=\xi _{3},\;\;\dot{\xi}_{3}=\mu ,
\label{chain}
\end{equation}
where $\mu \in {\mathbb R}$ is the control input. Let
$c_{1},c_{2}>0$ be given and
\begin{equation}
\delta \leq \min \left\{ \frac{4}{3},c_{1}\right\} .  \label{delta cond}
\end{equation}
Define the change of variables
\begin{equation}
z_{1}=\delta \left( c_{1}\xi _{1}+\xi _{2}\right) +\frac{\delta }{c_{1}}\xi
_{3},\,\,z_{2}=c_{1}\xi _{2}+\xi _{3},\,\,z_{3}=\xi _{3}  \label{z}
\end{equation}
and the functions
\begin{eqnarray}
U(z_{2},z_{3}) &=&\frac{1}{2}z_{3}^{2}+\int_{0}^{z_{2}}\sigma (s)ds
\label{UU} \\[0.4em]
g(z_{1},z_{2},z_{3}) &=&\frac{\delta }{c_{1}}\sigma (z_{1})+4\left(
1+U(z_{2},z_{3})\right)  \nonumber \\
&&\times \left( z_{3}+\sigma (z_{2})\right) -\frac{\delta }{2}\left(
z_{2}+z_{3}\right) .  \label{g}
\end{eqnarray}
Then the function
\begin{eqnarray}
V(\xi _{1},\xi _{2},\xi _{3}) &=&4U(z_{2},z_{3})+2U^{2}(z_{2},z_{3})
\nonumber \\
&&-\frac{\delta }{2}z_{2}z_{3}+\int_{0}^{z_{1}}\sigma (s)ds  \label{CLF}
\end{eqnarray}
is a CLF for (\ref{chain}) since its derivative along the trajectories of
system (\ref{chain}) in closed loop with the control
\begin{equation}
\mu (\xi _{1},\xi _{2},\xi _{3})=-c_{1}\sigma (z_{2})-c_{2}\sigma
(g(z_{1},z_{2},z_{3}))  \label{ubar}
\end{equation}
is negative definite. In particular, the bounded control (\ref{ubar})
renders (\ref{chain}) GAS to the origin.
\end{lemma}

\begin{remark}
\rm%
The saturation function $\sigma $ from Lemma \ref{mainlemma} is Lipschitz
continuous but not continuously differentiable. We chose this hard
saturation for simplicity of analysis. However, one can prove a similar
result with smoother saturation functions. For example, if one instead
chooses the $\mathcal{C}^{1}$ function
\[
\sigma (s)=\left\{
\begin{array}{l}
s\quad \mathrm{if}\;|s|\le 1 \\[-.8em]
\\
\mathrm{sign}(s)\left[ 1+\frac{1}{\pi }\tan ^{-1}\left( \pi \left(
\left| s\right| -1\right) \right) \right] \quad \mathrm{if}\;|s|>1,
\end{array}
\right.  \label{sigma} \]then the only change in the statement of
Lemma \ref{mainlemma} is to the right-hand side of the inequality
(\ref{delta cond}). The proof for this saturation is similar to the
proof we give below.
\end{remark}

\begin{remark}
\rm%
Condition (\ref{delta cond}) in Lemma \ref{mainlemma} allows cases where $%
c_{1}$ is arbitrarily large and $\delta $ is as small as desired, as well as
cases where $\max \{c_{1},c_{2},\delta \}$ is arbitrarily small, including
the case where $c_{1}=c_{2}=\delta \in (0,4/3]$.
\end{remark}

\begin{proof}%
We use (\ref{z}) and $\mu =-c_{1}\sigma (z_{2})+\alpha $, where $\alpha $ is
the new control input, to transform (\ref{chain}) into
\begin{equation}
\begin{array}{l}
\dot{z}_{1}=\delta \left( z_{2}-\sigma (z_{2})\right) +\frac{\delta }{c_{1}}%
\alpha \\
\dot{z}_{2}=c_{1}z_{3}-c_{1}\sigma (z_{2})+\alpha \\
\dot{z}_{3}=-c_{1}\sigma (z_{2})+\alpha .
\end{array}
\label{zdot}
\end{equation}

We first perform a Lyapunov-type analysis on the $(z_{2},z_{3})$-subsystem
using the function
\begin{equation}
V_{1}(z_{2},z_{3})=4U(z_{2},z_{3})+2U^{2}(z_{2},z_{3})-\frac{\delta }{2}%
z_{2}z_{3},  \label{V1}
\end{equation}
where $U(z_{2},z_{3})$ was defined in (\ref{UU}). By separately considering
the cases where $|z_{2}|>1$ and $|z_{2}|\leq 1$, one can easily verify that (%
\ref{V1}) is positive-definite and radially unbounded for $\delta <2$. In
fact, if $|z_{2}|\leq 1$, then
\begin{eqnarray}
V_{1}(z_{2},z_{3}) &\geq &4U(z_{2},z_{3})-\left| z_{2}z_{3}\right| \geq
2z_{3}^{2}+2z_{2}^{2}-\left| z_{2}z_{3}\right|  \nonumber \\
&\geq &\frac{3}{2}\left( z_{2}^{2}+z_{3}^{2}\right) .  \label{bound1 V1}
\end{eqnarray}
On the other hand, if $|z_{2}|>1$, then $\int\nolimits_{0}^{z_{2}}\sigma
(s)ds=1/2+\left( \left| z_{2}\right| -1\right) \geq \left| z_{2}\right| /2$,
so
\begin{eqnarray}
V_{1}(z_{2},z_{3}) &\geq &4\left( \frac{1}{2}z_{3}^{2}+\frac{1}{2}\left|
z_{2}\right| \right) +2\left( \frac{1}{2}z_{3}^{2}+\frac{1}{2}\left|
z_{2}\right| \right) ^{2}  \nonumber \\
&&-\left| z_{2}z_{3}\right|  \nonumber \\
&\geq &2z_{3}^{2}+2\left| z_{2}\right| +\frac{1}{2}z_{2}^{2}-\left|
z_{2}z_{3}\right| \geq \frac{3}{2}z_{3}^{2}  \nonumber \\
&&+2\left| z_{2}\right| .  \label{bound2 V1}
\end{eqnarray}

The time derivative of (\ref{V1}) along the trajectories of the $%
(z_{2},z_{3})$-subsystem is given by
\begin{eqnarray}
\dot{V}_{1} &=&4\left( 1+U(z_{2},z_{3})\right) \left( -c_{1}\sigma
^{2}(z_{2})+\left( z_{3}+\sigma (z_{2})\right) \alpha \right)  \nonumber \\
&&-\frac{\delta }{2}\left( c_{1}z_{3}-c_{1}\sigma (z_{2})+\alpha \right)
z_{3}  \nonumber \\
&&-\frac{\delta }{2}z_{2}\left( -c_{1}\sigma (z_{2})+\alpha \right)
\nonumber \\
&=&-4c_{1}\left( 1+U(z_{2},z_{3})\right) \sigma ^{2}(z_{2})-\frac{\delta
c_{1}}{2}z_{3}^{2}  \nonumber \\
&&+\frac{\delta c_{1}}{2}\sigma (z_{2})z_{3}+\frac{\delta c_{1}}{2}%
z_{2}\sigma (z_{2})  \nonumber \\
&&+\left[ 4\left( 1+U(z_{2},z_{3})\right) \left( z_{3}+\sigma (z_{2})\right)
\right.  \nonumber \\
&&\left. -\frac{\delta }{2}z_{3}-\frac{\delta }{2}z_{2}\right] \alpha .
\label{V1dot1}
\end{eqnarray}
Since $|\sigma (z_{2})|\leq |z_{2}|$ and $\left| z_{2}\sigma (z_{2})\right|
= $ $z_{2}\sigma (z_{2})$ give
\begin{equation}
-\frac{1}{4}z_{3}^{2}+\frac{1}{2}\sigma (z_{2})z_{3}\leq \frac{1}{4}\sigma
^{2}(z_{2})\le \frac{1}{4}z_{2}\sigma (z_{2}),  \label{ineq}
\end{equation}
(\ref{V1dot1}) becomes
\begin{eqnarray}
\dot{V}_{1} &\leq &-4c_{1}\left( 1+U(z_{2},z_{3})\right) \sigma ^{2}(z_{2})-%
\frac{\delta c_{1}}{4}z_{3}^{2}  \nonumber \\
&&+\frac{3}{4}\delta c_{1}z_{2}\sigma (z_{2})  \nonumber \\
&&+\left[ 4\left( 1+U(z_{2},z_{3})\right) \left( z_{3}+\sigma (z_{2})\right)
\right.  \nonumber \\
&&\left. \;\;-\frac{\delta }{2}z_{3}-\frac{\delta }{2}z_{2}\right]
\alpha . \label{V1dot2}
\end{eqnarray}
We now distinguish between two cases:

\begin{enumerate}
\item  \noindent $|z_{2}|\leq 1$. Then $\sigma (z_{2})=z_{2}$, $z_{2}\sigma
(z_{2})=\sigma ^{2}(z_{2})$, and since $\delta \leq 4/3$, we get
\begin{eqnarray}
\dot{V}_{1} &\leq &-3c_{1}\left( 1+U(z_{2},z_{3})\right) \sigma ^{2}(z_{2})-%
\frac{\delta c_{1}}{4}z_{3}^{2}  \nonumber \\
&&+\left[ 4\left( 1+U(z_{2},z_{3})\right) \left( z_{3}+\sigma (z_{2})\right)
\right.  \nonumber \\
&&\left. \;\;-\frac{\delta }{2}z_{3}-\frac{\delta }{2}z_{2}\right]
\alpha . \label{V1dot3}
\end{eqnarray}

\item  $|z_{2}|>1$. Then $\sigma (z_{2})=$ sign$(z_{2})$, $z_{2}\sigma
(z_{2})=\left| z_{2}\right| $, and
\begin{equation}
\renewcommand{\arraystretch}{1.4}
\begin{array}{l}
\left( 1+U(z_{2},z_{3})\right) \sigma ^{2}(z_{2})\;\ge
\;1+\int_{0}^{z_{2}}\sigma (s)ds\nonumber \\
\geq 1+\int_{0}^{\mathrm{sign}(z_{2})}sds+\int_{\mathrm{sign}%
(z_{2})}^{z_{2}}\sigma (s)ds\nonumber \\
=1+\int_{0}^{\mathrm{sign}(z_{2})}sds+\mathrm{sign}(z_{2})\int_{\mathrm{sign}%
(z_{2})}^{z_{2}}ds\nonumber \\
=\frac{1}{2}+|z_{2}|.
\end{array}
\label{integral}
\end{equation}

Since $\delta \leq 4/3$, applying (\ref{integral}) to (\ref{V1dot2}) again
yields (\ref{V1dot3}). Thus (\ref{V1dot3}) holds for all $z_{2}$.
\end{enumerate}

Finally, we perform a Lyapunov-type analysis on the whole system (\ref{zdot}%
) using the positive definite radially unbounded function $V$ from (\ref{CLF}%
). Since $V(z_{1},z_{2},z_{3})=V_{1}(z_{2},z_{3})+\int_{0}^{z_{1}}\sigma
(s)ds$, where $V_{1}(z_{2},z_{3})$ was defined in (\ref{V1}), it follows
from (\ref{V1dot3}) that the time derivative of (\ref{CLF}) along the system
trajectories satisfies
\begin{eqnarray}
\dot{V} &\leq &-3c_{1}\left( 1+U(z_{2},z_{3})\right) \sigma ^{2}(z_{2})-%
\frac{\delta c_{1}}{4}z_{3}^{2}  \nonumber \\
&&+\delta |z_{2}-\sigma (z_{2})|  \nonumber \\
&&+\left[ \frac{\delta }{c_{1}}\sigma (z_{1})+4\left(
1+U(z_{2},z_{3})\right) \left( z_{3}+\sigma (z_{2})\right) \right.  \nonumber
\\
&&\left. -\frac{\delta }{2}z_{3}-\frac{\delta }{2}z_{2}\right] \alpha .
\label{Vdot1}
\end{eqnarray}
We again distinguish between two cases:

\begin{enumerate}
\item  \noindent $|z_{2}|\leq 1$. Then $z_{2}-\sigma (z_{2})=0$, so (\ref
{Vdot1}) gives
\begin{eqnarray}
\dot{V} &\leq &-2c_{1}\left( 1+U(z_{2},z_{3})\right) \sigma ^{2}(z_{2})-%
\frac{\delta c_{1}}{4}z_{3}^{2}  \nonumber \\
&&+\left[ \frac{\delta }{c_{1}}\sigma (z_{1})+4\left(
1+U(z_{2},z_{3})\right) \right.  \nonumber \\
&&\left. \times \left( z_{3}+\sigma (z_{2})\right) -\frac{\delta }{2}z_{3}-%
\frac{\delta }{2}z_{2}\right] \alpha .  \label{Vdot2}
\end{eqnarray}

\item  $|z_{2}|>1$. Then (\ref{integral}) holds and since $|z_{2}-\sigma
(z_{2})|\leq |z_{2}|$ and $\delta \le c_{1}$, we again get (\ref{Vdot2}).
\end{enumerate}

\noindent Thus, (\ref{Vdot2}) holds for all $z_{2}$. Using $\alpha
=-c_{2}\sigma (g(z_{1},z_{2},z_{3}))$ and (\ref{g}), (\ref{Vdot2}) becomes
\begin{eqnarray}
\dot{V} &\leq &-2c_{1}\left( 1+U(z_{2},z_{3})\right) \sigma ^{2}(z_{2})-%
\frac{\delta c_{1}}{4}z_{3}^{2}  \nonumber \\
&&-c_{2}\sigma (g(z_{1},z_{2},z_{3}))g(z_{1},z_{2},z_{3})  \nonumber \\
&\leq &-2c_{1}\sigma ^{2}(z_{2})-\frac{\delta c_{1}}{4}z_{3}^{2}  \nonumber
\\
&&-c_{2}\sigma ^{2}(g(z_{1},z_{2},z_{3}))  \label{Vdot5}
\end{eqnarray}
by separately considering the cases $\left| g\right| $ $\geq 1$ and $\left|
g\right| <1$, i.e., $\dot{V}$ is negative definite. One can easily verify
that the origin is a unique equilibrium point of (\ref{zdot}) in closed loop
with $\alpha =-c_{2}\sigma (g(z_{1},z_{2},z_{3}))$. Thus, $%
(z_{1},z_{2},z_{3})=0$ is GAS \cite{khalil}. It follows from (\ref{z}) that $%
(\xi _{1},\xi _{2},\xi _{3})=0$ is GAS, and that (\ref{CLF}) is a CLF for (%
\ref{chain}). \end{proof}

\section{Main Result}

\label{sec4} We are now ready to state our main result on the global
asymptotic stabilization of (\ref{xdot amb}) with control saturation.

\begin{theorem}
\rm%
\label{mainthm} The control law
\begin{equation}
u(x)=\frac{\mu (x_{1},x_{2},X_{3})}{\beta _{0}+2\beta _{1}|x_{3}|},
\label{u}
\end{equation}
where $\mu (\cdot )$ was defined in (\ref{ubar}) and
\begin{equation}
X_{3}=\beta _{0}x_{3}+\beta _{1}x_{3}|x_{3}|,  \label{X3}
\end{equation}
in closed loop with (\ref{xdot amb}) ensures $x=0$ is GAS, and has
arbitrarily small amplitude. Moreover, with the choice (\ref{CLF}), the
function $x\mapsto V(x_{1},x_{2},X_{3})$ is a CLF for (\ref{xdot amb}).
\end{theorem}

\begin{proof}%
We use the change of variable (\ref{X3}) and the change of input (\ref{u})
to transform (\ref{xdot amb}) into the chain of integrators
\begin{equation}
\dot{x}_{1}=x_{2},\quad \dot{x}_{2}=X_{3},\quad \dot{X}_{3}=\mu .
\label{chain int}
\end{equation}
After setting $\mu $ to (\ref{ubar}) with $\xi _{1}=x_{1}$, $\xi _{2}=x_{2}$%
, and $\xi _{3}=X_{3}$, we can invoke Lemma \ref{mainlemma} to conclude $%
(x_{1},x_{2},X_{3})=0$ is GAS for (\ref{chain int}). It then follows from (%
\ref{X3}) that $x=0$ is GAS for the system (\ref{xdot amb}) in closed loop
with the feedback (\ref{u})-(\ref{X3}). The fact that $x\mapsto $ $%
V(x_{1},x_{2},X_{3})$ is a CLF for the AMB dynamics follows from the
proof of Lemma \ref{mainlemma}. \end{proof}

\begin{remark}
\rm%
Notice that the control law (\ref{u}) has arbitrarily small amplitude since
along the closed-loop trajectories,
\begin{equation}
\left| u(x(t))\right| \leq \frac{\left| \mu (x(t))\right| }{\beta _{0}}\leq
\frac{c_{1}+c_{2}}{\beta _{0}}\quad \forall t\geq 0.  \label{u ampl}
\end{equation}
Thus, one can choose $c_{1},c_{2}$ such that
\begin{equation}
c_{1}+c_{2}\leq \frac{\beta _{0}v_{\max }}{N}  \label{c12 cond}
\end{equation}
to ensure the AMB voltages satisfy the voltage constraint $\left|
v_{i}(t)\right| \leq v_{\max }$, $i=1,2$ $\forall t\geq 0$. From (\ref{c12
cond}), one can also see that the control gains can be adjusted to
accommodate sufficiently small bias levels due to the direct dependency of $%
\beta _{0}$ on the bias flux $\Phi _{0}$ (i.e., a low-bias control).
However, one cannot set $\Phi _{0}=0$ since the control may have a
singularity.
\end{remark}

\begin{remark}
\rm%
It is worth comparing our construction from Theorem \ref{mainthm} with the
AMB controller from Section VI of \cite{TA02}, where the stabilizing
controller took the form of nested hard saturations. In \cite{TA02}, the
proof that the control stabilizes the AMB system is based on the method of
\cite{teel1}, which uses the hardness of the saturation $\sigma $ in an
essential way. The argument shows that the nested hard saturations cause the
AMB dynamics (\ref{xdot amb}) to assume a simplified form after sufficiently
large time. Then one shows that this new system of `asymptotic' equations is
GAS using a Lyapunov-like analysis. Unfortunately, the Lyapunov function
used in the analysis of \cite{TA02} is not a CLF for the AMB\ dynamics (\ref
{xdot amb}). On the other hand, our construction is not restricted to the
standard hard saturation as explained in Remark 1, and thus can lead to
smoother stabilizing controllers. Also, the Lyapunov function from Theorem
\ref{mainthm} is negative definite along the closed-loop AMB trajectories
with the feedback (\ref{u}) \textit{for all time}, and thus is a CLF for the
AMB dynamics (\ref{xdot amb}).
\end{remark}

\begin{remark}
\rm%
Knowing a global CLF $V(x)$ and corresponding globally stabilizing feedback $%
u(x)$ for the AMB dynamics, we can construct a feedback $K(x)$ that renders
the system input-to-state stable (ISS) to exogenous disturbances \cite
{sontag1}. Choosing $V(x)$ and $u(x)$ from our theorem, and letting $\dot{x}%
=F(x,u)$ denote the open-loop AMB dynamics, it follows from [9, Theorem 1]
that the system $\dot{x}=$ $F(x,K(x)+d)$ is ISS relative to relative to any
measurable locally essentially bounded (e.g., locally bounded piecewise
continuous), exogenous disturbance $d(t)$ under the combined feedback $%
K(x)=u(x)-\partial V(x)/\partial x_{3}$. This is because the term $\partial
V(x)/\partial x_{3}$ in the feedback is the Lie derivative of $V$ in the
direction of the constant vector field $x\longmapsto (0,0,1)^{\intercal }$
that multiplies the input in the AMB dynamics. Note that in this case the
feedback $K(x)$ is not necessarily saturated. This is because if $K(x)$ were
any saturated feedback rendering $\dot{x}=F(x,d)$ ISS, and if $\bar{K}$ were
a positive global bound on $K(x)$, then the third component of the dynamic $%
\dot{x}_{3}=K(x)+d$ would have unbounded trajectories when the disturbance
is for example $d(t)\equiv 2\bar{K}$ , thus violating the ISS requirement.
Hence, there is no globally saturated feedback that renders the AMB dynamics
ISS to exogenous disturbances. Moreover, the CLF allows us to construct a
feedback $\tilde{K}(x)$ rendering the noisy system $\dot{x}=$ $F(x,\tilde{K}%
(x+e)+d)$ ISS with respect to the disturbance $d$ for small measurement
errors $e(t)$ in the controller \cite{sontag2}. In the context of AMBs, $d$
could represent rotor unbalance or vibrations transmitted through the
bearing\ foundation while $e$ could represent inaccuracies in sensing
displacement, velocity, and/or flux.
\end{remark}

\section{Numerical Example\label{sec5}}

The proposed control law from Theorem \ref{mainthm} above was simulated on
the same high-fidelity model of the one degree-of-freedom magnetic
levitation system \cite{knospe,lindlau} used in \cite{TA02}. The model
includes flux leakage, magnetic material saturation, coil resistance,
voltage saturation, and flexible modes, and thus is a fairly accurate
representation of the actual system test rig \cite{lindlau}. We refer the
reader to \cite{lindlau} for a complete description of the system model.

As in \cite{TA02}, the system was simulated with $\Phi _{0}=10$ $\mu $Wb, $%
v_{\max }=10$ V, and initial conditions set to $y(0)=0.15$ mrad, $\dot{y}%
(0)=0$ mrad/s, $\Phi _{1}(0)=20$ $\mu $Wb, and $\Phi _{2}(0)=60$ $\mu $Wb.
The above flux values are for the \textit{total }flux in the magnetic
levitation system, i.e., including gap flux and flux leakage. The control
parameters $c_{1},c_{2},\delta $ were selected according to (\ref{c12 cond})
and (\ref{delta cond}) to yield the fastest convergence possible for the
state $y$. This resulted in the following parameter values: $c_{1}=5465$, $%
c_{2}=5300$, and $\delta =4/3$. The simulation results for the system states
$y$, $\dot{y}$, and $\Phi _{gi}$ ($i=1,2$), and for the control voltages $%
v_{i}$ ($i=1,2$) are shown in Figure 2 below, where $\Phi _{gi}$
denotes the gap flux. When compared with the simulation of the
controller from Section VI of \cite{TA02} (see Figure 8 of
\cite{TA02}), the our control yields similar closed-loop performance
although it can be tuned to produce a faster response as seen from
Figure 2.

\section{Conclusion}

\label{sec6}We constructed a stabilizing controller and control-Lyapunov
function for low-bias active magnetic bearings (AMBs) with voltage
saturation. The AMB\ dynamics used in the construction is borrowed from a
recently-developed model for AMBs operating with only one electromagnet
active at any given time. Our construction is based on a new forwarding-like
approach to stabilizing a three-dimensional chain of integrators with input
saturation. The proposed AMB control law has arbitrarily small amplitude and
the flexibility of being implemented with different saturation-like
functions.

\begin{figure}[ht]
\begin{center}
\scalebox{.175}{\includegraphics{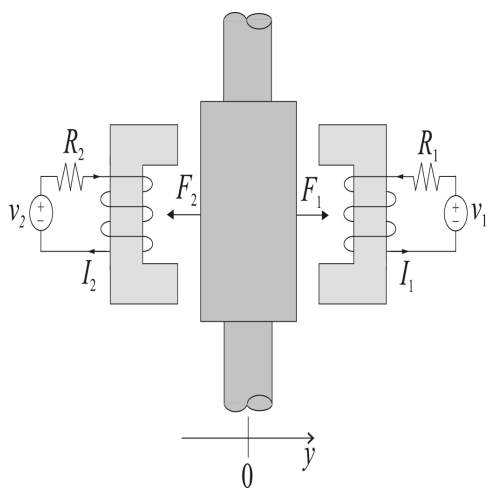}}
\end{center}
\caption{ Schematic of the AMB system.}
\end{figure}

\begin{figure}[ht]
\begin{center}
\scalebox{.52}{\includegraphics{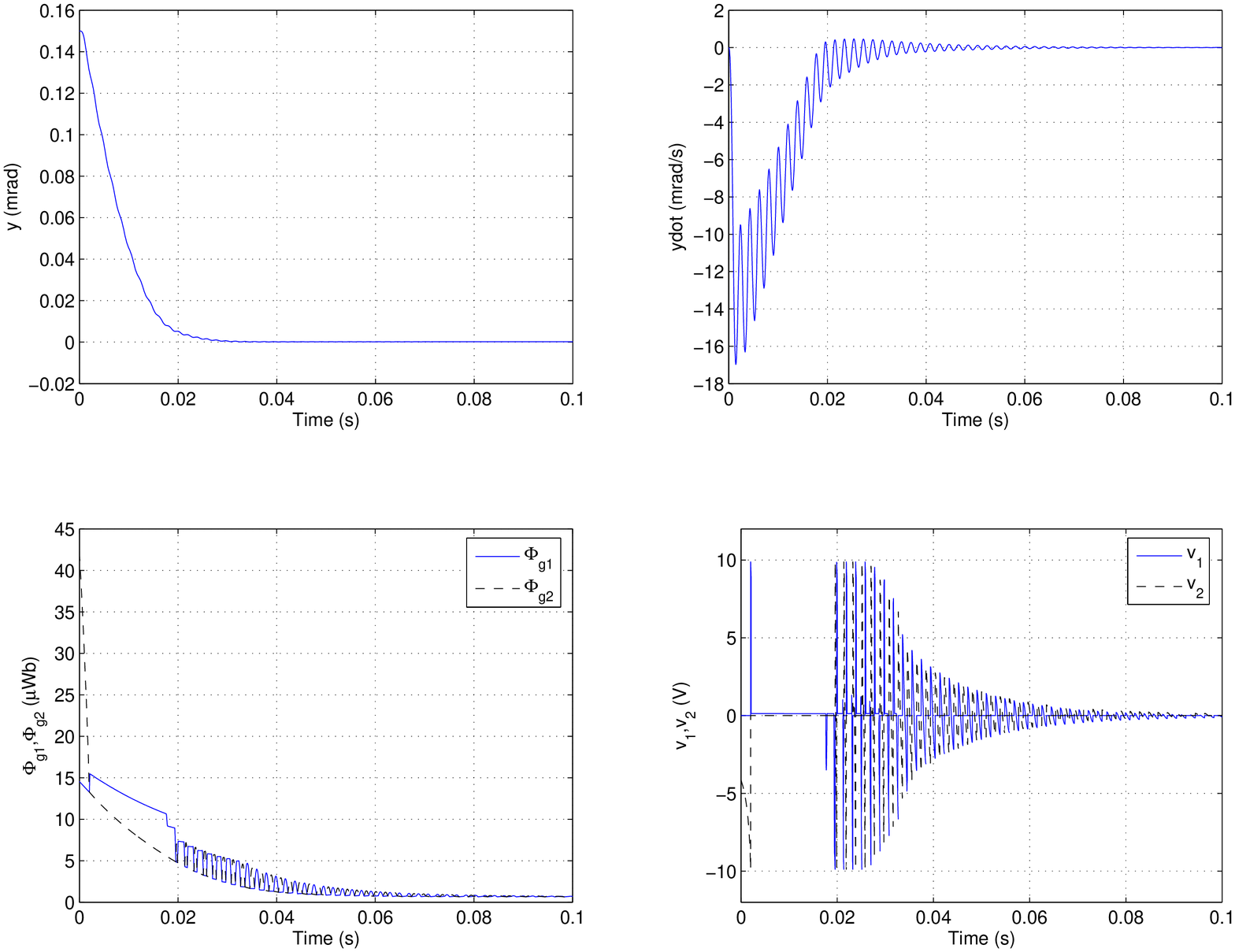}}
\end{center}
\caption{State trajectories and control voltages.}
\end{figure}

\begin{thebibliography}{99}
\bibitem{fujita}  M. Fujita, K. Hatake, and F. Matsumura, ``Loop Shaping
Based Robust Control of a Magnetic Bearing,'' \textit{IEEE Control Syst. Mag.%
}, Vol. 13, No. 4, pp. 57-65, 1993.

\bibitem{hu}  T. Hu, Z. Lin, and P.E. Allaire, ``Power-Loss Minimization in
Magnetic Bearing Systems,'' \textit{Proc. IEEE Conf. Decision and Control},
pp. 854-859, Las Vegas, NV, 2002.

\bibitem{khalil}  H. Khalil, \textit{Nonlinear Systems}, Upper Saddle River,
NJ: Prentice Hall, 2002.

\bibitem{knospe}  C. Knospe, ``The Nonlinear Control Benchmark System,''
Proc. American Control Conf., pp. 2134-2138, Chicago, IL, 2000.

\bibitem{krstic}  M. Krstic, I. Kanellakopoulos, and P. Kokotovic, \textit{%
Nonlinear and Adaptive Control Design}, New York, NY: John Wiley \& Sons,
1995.

\bibitem{lindlau}  J.D. Lindlau, \textit{Dynamic Force Biasing of Active
Magnetic Bearings via Feedback Linearization}, M.S. Thesis, University of
Virginia, 1999.

\bibitem{matsumura}  F. Matsumura and T. Yoshimoto, ``System Modeling and
Control Design of a Horizontal-Shaft Magnetic-Bearing System,'' \textit{IEEE
Trans. Magnetics}, Vol. MAG-22, No. 3, 1986.

\bibitem{mazenc}  F. Mazenc and L. Praly, ``Adding Integrations, Saturated
Controls, and Stabilization for Feedforward Systems,'' \textit{IEEE Trans.
Automatic Control}, Vol. 41, No.11, pp. 1559-1578, 1996.

\bibitem{mohamed1}  A.M. Mohamed and F.P. Emad, ``Conical Magnetic Bearings
with Radial and Thrust Control,'' \textit{IEEE Trans. Automatic Control},
Vol. 37, No. 12, pp. 1859-1868, 1992.

\bibitem{motee cdc03}  N. Motee and M.S. de Queiroz, ``A Switching Control
Strategy for Magnetic Bearings with a State-Dependent Bias,'' \textit{Proc.
IEEE Conf. Decision and Control}, pp. 245-250, Maui, HI, 2003.

\bibitem{sontag1}  E.D. Sontag, ``Input-to-State Stability: Basic Concepts
and Results,'' in \textit{Nonlinear and Optimal Control Theory},
Springer Lecture Notes in Mathematics, to appear,
http://www.math.rutgers.edu/ $\! \sim$sontag/papers.html.

\bibitem{sontag2}  E.D. Sontag, ``Clocks and Insensitivity to Small
Measurement Errors,'' \textit{ESAIM: Control, Optimisation, and Calculus of
Variations}, Vol. 4, pp. 537-557, 1999.

\bibitem{teel1}  A. Teel, ``Using Saturation to Stabilize a Class of
Single-Input Partially Linear Composite Systems,'' \textit{Proc. IFAC
Nonlinear Control Systems Design Symp.}, pp. 224-229, Bordeaux, France, 1992.

\bibitem{teel2}  A. Teel, ``A Nonlinear Small Gain Theorem for the Analysis
of Control Systems with Saturation,'' \textit{IEEE Trans. Automatic Control}%
, Vol. 41, No. 9, pp. 1256-1270, 1996.

\bibitem{tsiotras1}  P. Tsiotras and E. Velenis, ``Low-Bias Control of AMB's
Subject to Saturation Constraints,'' \textit{Proc. IEEE Conf. Control Appl.}%
, pp. 138-143, Ancorage, AK, 2000.

\bibitem{TA02}  P. Tsiotras and M. Arcak, ``Low-Bias Control of AMB Subject
to Voltage Saturation:\ State-Feedback and Observer Designs,'' \textit{IEEE
Trans. Control Syst. Tech.}, Vol. 13, No. 2, pp 262-273, 2005.

\bibitem{tsiotras4}  P. Tsiotras and B.C. Wilson, ``Zero and Low-Bias
Control Designs for Active Magnetic Bearings,'' \textit{IEEE Trans. Control
Syst. Tech.}, Vol. 11, No. 6, pp 889-904, 2003.

\bibitem{embook}  H.H. Woodson and J.R. Melcher, \textit{Electromechanical
Dynamics - Part I: Discrete Systems}, New York, NY: John Wiley \&
Sons, 1968.
\end{thebibliography}
\end{document}